\documentclass[11pt]{amsart}

\newtheorem{theorem}{Theorem}
\newtheorem{lem}[theorem]{Lemma}

\theoremstyle{definition}

\newtheorem{ex}[theorem]{Example}
\newtheorem{rem}[theorem]{Remark}

\newcommand{\eqr}[1]{~\mbox{$(${\rm \ref{#1}}$)$}}
\newcommand{\Q}{\mathbb{Q}}
\renewcommand{\i}{\mathbf{i}}
\renewcommand{\j}{\mathbf{j}}
\newcommand{\bt}{\mathbf{t}}
\newcommand{\s}{\mathbf{s}}
\newcommand{\y}{\mathbf{y}}
\newcommand{\bk}{\mathbf{k}}
\newcommand{\newdots}{\raisebox{-1mm}{.}.\raisebox{1mm}{.}}

\begin{document}

\title[Multiplicities of Points on Schubert
Varieties]{Multiplicities of Points on Schubert Varieties in
  Grassmannians}

\author{Joachim Rosenthal} \address{Department of Mathematics,
  University of Notre Dame, Notre Dame, Indiana 46556}
\email{Rosenthal.1@nd.edu} 
\thanks{The authors were partially
  supported by NSF grants DMS-9610389 and DMS-9625511.}

\author{Andrei Zelevinsky} 
\address{Northeastern University, Department of Mathematics,
  Boston, MA 02115} 
\email{andrei@neu.edu}

\date{January 14, 1999}

\subjclass{Primary: 14M15. Secondary: 14H20. }

\keywords{Schubert varieties, singularities, multiplicities,
  partial difference equation.}

\begin{abstract}
  We obtain an explicit determinantal formula for the
  multiplicity of any point on a classical Schubert variety.
\end{abstract}

\maketitle

\section{Main result}

An important invariant of a singular point on an algebraic
variety $X$ is its \emph{multiplicity}: the normalized leading
coefficient of the Hilbert polynomial of the local ring.  The
main result of the present note is an explicit determinantal
formula for the multiplicities of points on Schubert varieties in
Grassmannians.  This is a simplification of a formula obtained
in~\cite{ro86}.  More recently, the recurrence relations for
multiplicities of points on more general (partial) flag varieties
were obtained in~\cite{la95a,la90a}.  However, to the best of our
knowledge the case of Grassmannians remains the only case for
which an explicit formula for multiplicities is available.

Fix positive integers $d$ and $n$ with $0 \leq d \leq n$, and
consider the Grassmannian $Gr_d (V)$ of $d$-dimensional subspaces
in a $n$-dimensional vector space $V$ (over an algebraically
closed field of arbitrary characteristic).  Recall that Schubert
varieties in $Gr_d (V)$ are parameterized by the set $I_{d,n}$ of
integer vectors $\i = ( i_{1},\ldots ,i_{d})$ such that $1 \leq
i_{1} < \ldots <i_{d} \leq n$.  For a given complete flag $\{ 0\}
= V_0 \subset V_{1} \subset \cdots \subset V_{n}= V$, the
Schubert variety $X_\i$ is defined as follows:
$$
X_\i := \{W \in Gr_d(V)\mid \dim (W \bigcap V_{i_{k}}) \geq k
\ \mbox{ for } k=1,\ldots ,d\} \ .
$$
The Schubert cell $X_\i^0$ is an open subset in $X_\i$ given
by
$$
X_\i^0 := \{ W \in X_\i \mid \dim (W \bigcap V_{i_{k}-1}) =
k-1 \mbox{ for } k=1,\ldots ,d\} \ .$$
It is well known that the
Schubert variety $X_\i$ is the disjoint union of Schubert cells
$X_\j^0$ for all $\j\leq \i$ in the componentwise partial order
on $I_{d,n}$.  The multiplicity of a point $x$ in $X_\i$ is
constant on each Schubert cell $X_\j^0 \subset X_\i$, and we
denote this multiplicity by $M_\j(\i)$.

Our main result is the following explicit formula for $M_\j(\i)$
(where the binomial coefficients ${a \choose b}$ are subject to
the condition that ${a \choose b} = 0$ for $b < 0$):

\begin{theorem}
\label{main}
The multiplicity $M_\j(\i)$ of a point $x \in X_\j^0 \subset
X_\i$ is given by
  \begin{equation}
  \label{main-form}
    M_\j(\i)=(-1)^{s_1 + \cdots + s_d} \det \left[
  \begin{array}{cccc}
  {i_1\choose -s_1} & \ldots & \ldots& {i_d \choose -s_d}\\
  {i_1\choose 1-s_1} & \ldots & \ldots& {i_d\choose 1-s_d}\\
  \vdots &  & & \vdots\\
  {i_1\choose d-1-s_1} & \ldots & \ldots& {i_d\choose d-1-s_d}
  \end{array}
\right] \ ,
\end{equation}
where
\begin{equation}
  \label{sij}
s_q := \# \{ j_p \mid i_q < j_p\} \ .
\end{equation}
\end{theorem}

The proof of Theorem~\ref{main} will be given in the next
section.  Although determinants of matrices formed by binomial
coefficients were extensively studied by combinatorialists (see,
e.g., \cite{GV}), the experts whom we consulted did not recognize
the determinant in (\ref{main-form}).

We conclude this section by an example illustrating
Theorem~\ref{main}.

\begin{ex} 
\label{separated_ij}
Assume the indices $\i,\j$ satisfy $j_d \leq i_1$.  In this
situation the numbers $s_1,\ldots,s_d$ attain the smallest
possible value: $s_1=\cdots=s_d=0$.  Then the $(p,q)$-entry of
the determinant in (\ref{main-form}) has the form $P_p (i_q)$,
where $P_p (t)$ is a polynomial with the leading term
$t^{p-1}/(p-1)!$.  It follows that
  \begin{equation}
    M_\j(\i)= \frac{1}{1!\cdots (d-1)!}V(\i)= 
    \frac{1}{1!\cdots (d-1)!}\prod_{p>q}(i_p-i_q) \ ,
  \end{equation}
  where $V(\i)$ is the Vandermonde determinant $\det
  ((i_q^{p-1}))$.
\end{ex}

\section{Proof of Theorem \ref{main}}
Fix two vectors $\j\leq \i$ from $I_{d,n}$, and let
\begin{equation*}
  \deg (\j,\i):= d-\# \{ i_q \mid i_q \in \{ j_1,\ldots,j_d\}\} \ .
\end{equation*}
For a nonnegative integer vector $\s = (s_1, \ldots, s_d)$, we
set
\begin{equation*}
  |\s|:=s_1+\cdots+s_d \ .
\end{equation*}
As shown in~\cite{ro86} and~\cite[page 202]{la90a}, the
multiplicity $M_\j(\i)$ satisfies the initial condition
$M_\j(\j)=1$ and the partial difference equation
\begin{equation}
  \label{recu}
  M_\j(\i)=\frac{1}{\deg (\j,\i)}\sum_{\bk} M_\j(\bk) \ , 
\end{equation}
where the sum is over all $\bk \in I_{d,n}$ such that $\j\leq
\bk< \i$, and $|\bk|=|\i|-1$.

To prove\eqr{main-form}, we proceed by induction on $|\i|$.  The
initial step is to verify\eqr{main-form} for $\i = \j$.  In this
case the numbers $s_1,\ldots, s_d$ attain their maximum possible
value: $s_q=d-q$.  It follows that
\begin{equation}
  (-1)^{|\s|}\det\left[
  \begin{array}{cccc}
   0 & \ldots & 0 & 1\\
   \vdots &   & 1 & \ast\\
  0 &\newdots & \newdots &\vdots \vspace{3mm}\\
  1 & \ast & \ldots & \ast
  \end{array}
\right] =1 = M_\j(\j)\ ,
\end{equation}
as required.

For the inductive step, we introduce some notation.  To any
nonnegative integer vector $\s = (s_1, \ldots, s_d)$ we associate
a polynomial $P_\s (\bt) \in \Q[\bt] = \Q[t_1, \ldots, t_d]$
defined by
\begin{equation}
  \label{poly}
    P_\s (\bt)=(-1)^{|\s|} \det \left[
  \begin{array}{cccc}
  {t_1\choose -s_1} & \ldots & \ldots& {t_d \choose -s_d}\\
  {t_1\choose 1-s_1} & \ldots & \ldots& {t_d\choose 1-s_d}\\
  \vdots &  & & \vdots\\
  {t_1\choose d-1-s_1} & \ldots & \ldots& {t_d\choose d-1-s_d}
  \end{array}
\right] \ ;
\end{equation}
here ${t \choose s}$ is the polynomial $t(t-1) \cdots (t-s+1)/s!$
for $s \geq 0$, and ${t \choose s} = 0$ for $s < 0$.  Thus our
goal is to show that $M_\j (\i) = P_\s (\i)$ with $\s$ given by
(\ref{sij}).

For $q = 1, \ldots, d$, let $\Delta_q : \Q[\bt] \to \Q[\bt]$
denote the partial difference operator $\Delta_q P (\bt) = P(\bt)
- P(\bt - e_q)$, where $e_1, \ldots, e_d$ are the unit vectors in
$\Q^d$.  Here is the key lemma.

\begin{lem} 
\label{dif_Euler}
For any nonnegative integer vector $\s$, the corresponding
polynomial $P_\s (\bt)$ satisfies the partial difference equation
\begin{equation}
  \label{dif_eqn}
(\Delta_1 + \cdots + \Delta_d) P = 0 \ .
\end{equation}
\end{lem}

\begin{proof}
  First notice that the Vandermonde determinant $V(\bt) =
  \prod_{p>q}(t_p-t_q)$ satisfies (\ref{dif_eqn}) since it is a
  non-zero skew-symmetric polynomial of minimal possible degree,
  and the operator $\Delta_1 + \cdots + \Delta_d$ preserves the
  space of skew-symmetric polynomials.  The vector space of
  solutions of (\ref{dif_eqn}) is also invariant under
  translations $\bt \mapsto \bt + \bk$ so it is enough to show
  that each $P_\s (\bt)$ is a linear combination of polynomials
  $V(\bt + \bk)$.  Here is the desired expression:
\begin{equation}
\label{main-form2}
P_\s (\bt) = 
   \frac{1}{1!\cdots (d-1)!}\sum_{0 \leq \bk \leq \s} 
  (-1)^{|\bk|}{s_1\choose k_1}\cdots {s_d\choose k_d}V(\bt+\bk)\ .
  \end{equation}
  Let us prove (\ref{main-form2}).  The same argument as in
  Example~\ref{separated_ij} above shows that
\begin{equation}\label{V-form}
   \frac{1}{1!\cdots (d-1)!} V(\bt+\bk)=
\det \left[
  \begin{array}{cccc}
  {t_1+k_1\choose 0} & \ldots & \ldots& {t_d+k_d\choose 0}\\
  {t_1+k_1\choose 1} & \ldots & \ldots& {t_d+k_d\choose 1}\\
  \vdots &  & & \vdots\\
  {t_1+k_1\choose d-1} & \ldots & \ldots& {t_d+k_d\choose d-1}
  \end{array}
\right] \ .
  \end{equation}
  Substituting this expression into (\ref{main-form2}) and
  performing the multiple summation, we see that the right hand
  side becomes the determinant of the $d \times d$ matrix whose
  $(p,q)$-entry is
\begin{equation*}
\sum_{k_q=0}^{s_q}(-1)^{k_q}{s_q\choose k_q}{t_q+k_q \choose p-1}= 
(-1)^{s_q}{t_q\choose p-1-s_q}
\end{equation*}
(the last equality is a standard binomial identity).  This
completes the proof of (\ref{main-form2}) and
Lemma~\ref{dif_Euler}.
\end{proof}

One last piece of preparation before performing the inductive
step: the Pascal binomial identity ${t \choose s} = {t-1 \choose
  s} + {t-1 \choose s-1}$ implies that
\begin{equation}
\label{delta P}
\Delta_q P_\s (\bt) = - P_{\s + e_q} (\bt - e_q)
\end{equation}
for any nonnegative integer vector $\s$ and any $q = 1, \ldots, d$.

To conclude the proof of Theorem~\ref{main}, suppose that $\j <   \i$ and assume 
by induction that $M_\j (\bk)$ is given
  by\eqr{main-form} for any $\bk \in I_{d,n}$ such that $\j \leq \bk < \i$.
Let $\s$ be the vector given by (\ref{sij}). 
In view of \eqr{recu}, the desired equality $M_\j (\i) = P_\s (\i)$ 
is a consequence of the following:
\begin{equation}
\label{induction}
\deg (\j,\i) P_\s (\i) - \sum_{\bk} M_\j(\bk)  = 0 \ ,
\end{equation}
  where the sum is over all $\bk \in I_{d,n}$ such that
  $\j\leq \bk< \i$, and $|\bk|=|\i|-1$.

We shall deduce \eqr{induction} from the equality
  $$
  \sum_{q=1}^d \Delta_q P_\s (\i) = 0
  $$
  provided by Lemma~\ref{dif_Euler}.
To do this, we compute $\Delta_q P_\s (\i)$ in each of the 
following mutually exclusive cases (we use the
conventions $i_0 = 0$ and $s_0 = d$):

\noindent {\bf Case 1:} $i_q \notin \{j_1, \ldots, j_d\}$, 
  $i_q - 1 > i_{q-1}$.  
Then $\bk := \i - e_q$ belongs to $I_{d,n}$, and
we have $\j \leq \bk$.  
Replacing $\i$ by $\bk$ in (\ref{sij}) does not change the vector $\s$.  
By our inductive assumption, $P_\s (\bk) = M_\j (\bk)$, and so $\Delta_q
  P_\s (\i) = P_\s (\i) - M_\j (\bk)$.

\noindent {\bf Case 2:} $i_q \notin \{j_1, \ldots, j_d\}$, 
  $i_q - 1 = i_{q-1}$.  For such $q$, we have $P_\s (\i - e_q) =
  0$ since the corresponding determinant has the $(q-1)$th and
  $q$th columns equal to each other.  Thus $\Delta_q P_\s (\i) = P_\s (\i)$.

\noindent {\bf Case 3:} $i_q \in \{j_{q+1}, \ldots, j_d\}$, 
  $i_q - 1 > i_{q-1}$.  As in Case 1, we have $\bk :=
  \i - e_q \in I_{d,n}$, and $\j \leq \bk$.  However now
  replacing $\i$ by $\bk$ in (\ref{sij}) changes $\s$ to $\s +
  e_q$.  Combining the inductive assumption with (\ref{delta P}),
  we conclude that $\Delta_q P_\s (\i) = - P_{\s + e_q} (\bk) = -
  M_\j (\bk)$.

\noindent {\bf Case 4:} $i_q \in \{j_{q+1}, \ldots, j_d\}$, 
  $i_q - 1 = i_{q-1}$.
In this case, the $d \times d$ matrix whose determinant is 
$P_{\s + e_q}(\i - e_q)$ has the $(q-1)$th and $q$th columns equal to each
  other, hence $\Delta_q P_\s (\i) = - P_{\s + e_q} (\bk) = 0$.

\noindent {\bf Case 5:} $i_q = j_q$. 
Then we have
  $$
  s_1 \geq s_2 \geq\cdots\geq s_{q-1}\geq s_q + 1 = d+1-q \ ,
  $$
  and so the $d \times d$ matrix whose determinant is $P_{\s +
    e_q}(\i - e_q)$ has a zero $(d+1-q) \times q$ submatrix.
As in Case 4, this implies $\Delta_q P_\s (\i) = - P_{\s + e_q} (\bk) = 0$.

Adding up the contributions $\Delta_q P_\s (\i)$ from all these cases,
we obtain \eqr{induction}; this completes the proof of Theorem~\ref{main}.

\begin{rem}
  In~\cite{ro86}, the multiplicity $M_\j (\i)$ was expressed as a
  multiple sum given by\eqr{main-form2}.
\end{rem}

\begin{rem}
  The multiplicity $M_\j (\i)$ is by definition a positive integer.
  The partial difference equation \eqr{recu} (combined with the initial condition 
  $M_\j (\j) = 1$) makes the positivity of $M_\j (\i)$ obvious but the fact 
  that $M_\j (\i)$ is an integer becomes rather mysterious.
  On the other hand, Theorem~\ref{main} makes it clear that 
  $M_\j (\i)$ is an integer but not that $M_\j (\i) > 0$. 
  It would be interesting to find an expression for $M_\j (\i)$
  that makes obvious both properties.  
\end{rem}

\begin{rem}
  The space of all polynomial solutions of the partial difference
  equation\eqr{dif_eqn} can be described as follows.  Let $\y =
  (y_1, \ldots, y_d)$ be an auxiliary set of variables, and let
  $\varphi: \Q[\y] \to \Q[\bt]$ be the isomorphism of vectors
  spaces that sends each monomial $\prod_{q=1}^d y_q^{n_q}$ to
  $\prod_{q=1}^d t_q (t_q + 1) \cdots (t_q + n_q -1)$.  The map
  $\varphi$ intertwines each $\Delta_q$ with the partial
  derivative $\frac{\partial}{\partial y_q}$.  It follows that
  the space of solutions of\eqr{dif_eqn} is the image under
  $\varphi$ of the $\Q$-subalgebra in $\Q[\y]$ generated by all
  differences $y_p - y_q$.
\end{rem}

\begin{rem}
  Jerzy Weyman informed us about the following determinantal
  formula (unpublished) for the multiplicity $M_\j (\i)$ in the
  special case when $\j = (1, 2, \ldots, d)$.  Let $\lambda$ be
  the partition $(i_d - d, \ldots, i_2 - 2, i_1 - 1)$, and let
  $\lambda = (\alpha_1, \ldots, \alpha_r|\beta_1, \ldots,
  \beta_r)$ be the Frobenius notation of $\lambda$
  (see~\cite{mac}).  According to J.~Weyman, $M_\j (\i)$ is equal
  to the determinant of the $r \times r$ matrix whose
  $(p,q)$-entry is ${\alpha_p + \beta_q \choose \alpha_p}$.  It
  is not immediately clear why this determinantal expression
  agrees with the one given by\eqr{main-form}.
\end{rem}

\section*{Acknowledgements}

We are grateful to V.~Lakshmibai who initiated this project by
suggesting to one of us (J.~R.) to publish the results of his
thesis~\cite{ro86}. We thank Sergey Fomin, Ira Gessel and Jerzy
Weyman for helpful conversations.


\providecommand{\bysame}{\leavevmode\hbox
  to3em{\hrulefill}\thinspace}

\end{document}